\documentclass{elsarticle}

\newcommand{\be}{\begin{eqnarray}}
\newcommand{\ee}{\end{eqnarray}}

\newcommand{\short}{\mbox{Short}}
\newcommand{\cut}{\mbox{Cut}}

\newtheorem{theorem}{Theorem}

\newtheorem{lemma}{Lemma}

\newcommand{\QED}{\nobreak \ifvmode \relax \else
      \ifdim\lastskip<1.5em \hskip-\lastskip
      \hskip1.5em plus0em minus0.5em \fi \nobreak
      \vrule height0.75em width0.5em depth0.25em\fi}
	  
\newcommand{\myfig}[4]{\centerline{\resizebox{!}{#1\textwidth}{\includegraphics{#2}}}
\caption{#4}\label{#3}}
	  
\usepackage{graphicx}
\usepackage{amssymb}

\begin{document}

\title{Minimal and irreducible links in the Shannon game}

\author{Andrew M. Steane}
\address{
Centre for Quantum Computing,
Department of Atomic and Laser Physics, Clarendon Laboratory,\\
Parks Road, Oxford, OX1 3PU, England.
}

\date{\today}

\begin{abstract}
We discuss weak and strong links (`virtual connections') in the Shannon game.
General properties of these links are discussed, leading to
a method to find all links of given size by a suitably pruned exhaustive search.
This is applied to links on graphs of up to 11 vertices. We discuss
the concept of reducibility of such links. Three simple reductions
are considered, including one, the `short-cut', not previously described.
The complete sets of irreducible weak links on up to 11 vertices 
and strong links on up to 10 vertices are
presented. Some applications to the analysis of Hex are noted.
\end{abstract}

\begin{keyword}
Shannon game  \sep weak link \sep Hex \sep virtual connection \sep game theory
\end{keyword}

\maketitle

\section{Introduction}

The {\em Shannon game} is a vertex-colouring game played on a simple graph.
It leads to a variety of interesting mathematical problems and concepts;
see \cite{06Hayward,Arneson:2010:SHB:1950322.1950323} for reviews. There
is always a winner (the game cannot
be drawn), but the problem of determining the winner of a general game 
is PSPACE-complete. In the absence of an efficient algorithm for this
general task, the most natural strategy is to try to break down the task
into sub-games. This leads to the concept of `strong' and `weak' 
{\em virtual connections} or {\em links}, and to the concept of
the {\em multi-Shannon game}, see \cite{06Hayward,06Bj,12steane}. 
A `link' (we prefer this term to `virtual connection') is a 
set of vertices and edges which suffices to allow one player (\short)
to guarantee to form a path between two vertices at the
ends of the link, following the rules of the game. It is essentially
another Shannon game. If \short\ has a 2nd-player winning strategy
for such a link, it is said to be a {\em strong link}; if \short\ has a 1st-player
but not a 2nd-player winning strategy, it is said to be a {\em weak link}.
Figure \ref{f.minweak} shows a collection of small weak links. We will
refer to some of these by the names W1, W3, W5X, etc. as shown in the figure.
W1 is the smallest weak link.


\begin{figure}
\myfig{0.6}{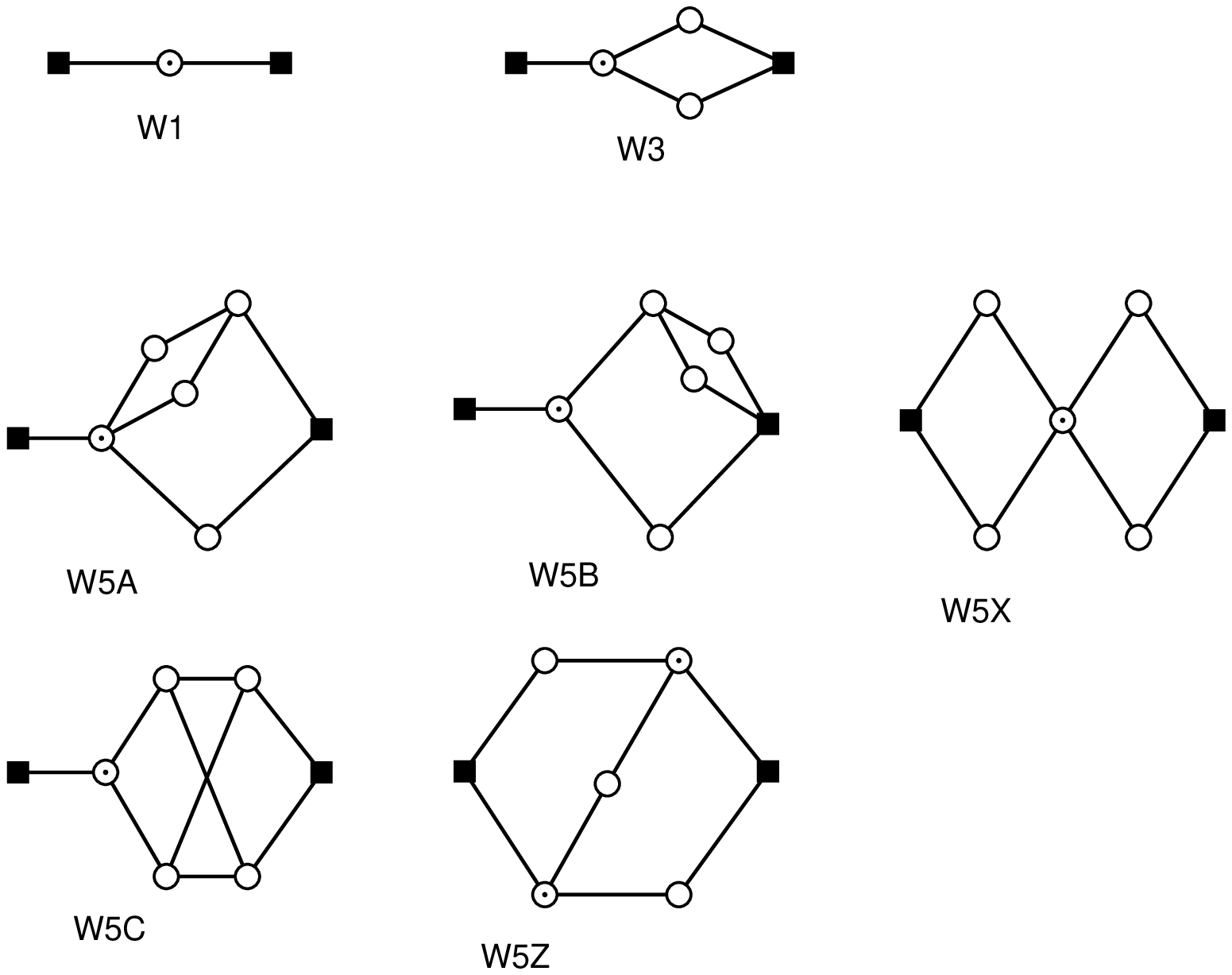}{f.minweak}
{The set of all minimal weak links of weight up to 5, with names
to allow them to be easily referred to in the text. A dot indicates a pivot.}
\end{figure}

The {\em weight} of a link is the number of vertices in it, excluding
the terminals. W1 has weight 1. The smallest strong link has zero weight;
it is a direct adjacency between two vertices. The next smallest strong link
is the following, called a `bridge' or S2:

\begin{picture}(40,30)(130,0)
\multiput(147.5,7.5)(60,0){2}{\rule{5pt}{5pt}}
\multiput(152.5,11)(30,-10){2}{\line(3,1){25}}
\multiput(152.5,9)(30,10){2}{\line(3,-1){25}}
\multiput(180,19)(0,-18){2}{\circle{5}}
\end{picture}

One strategy to solve the Shannon game is to first identify all small links
in a given game, and then use these to identify larger links, until eventually
a link reaching all the way from one terminal of the game to the other is found,
or else the non-existence of such a link is established. In practice this strategy
is not normally used on its own, because the task of identifying large links
is prohibitively difficult, but it is an important part of work on the game.

Early exponents of the game noticed the `bridge' and also the `chain' which
is a sequence of bridges in series, and various other combinations. The simplest
combinations (`OR rule', `AND rule') were described by Anshelevich \cite{00Anshelevich}
and implemented in a computer program {\em Hexy}. The `OR rule' is the statement that
a strong link exists between given vertices 
whenever there is a set of weak links between those vertices having empty
total intersection (that is, for every vertex there is at least one weak link
not involving that vertex). It is easy to see that this is both a 
sufficient and a necessary
condition: if a strong link exists then, after \cut\ has made a move at some
vertex $u$, \short\ now has a 1st-player winning strategy on the remaining vertices,
hence a weak link that does not involve $u$, and this must be true for all $u$
in the strong link. It follows that {\em to find all strong and weak links
in the Shannon game, it suffices to find only weak links and apply the OR-rule}. 

Work to date has been of two broad types. One is to investigate various strategies
and concepts, not all of them based on weak and strong links, see for example
\cite{00Anshelevich,01Yang,03Yang,06Hayward}. The other is
to use a combination of strategies to solve larger and larger games,
especially Hex problems, by computer \cite{Arneson:2010:SHB:1950322.1950323}.
In the course of all this work, weak links have been discussed and their simplest
examples found, but the general problem of an exhaustive analysis has
not been attempted. 

In this paper we investigate the structure of links in general.
We first ask the question: what weak links
exist? In other words, we would like to know about the complete set of all
weak links of given size. To focus attention on the essential structure of
these links, we consider the {\em minimal} links---those which have
no extraneous edges (a formal definition is given below). 
We present various general properties of such links
which help to reduce an exhaustive search for them. We have carried out
such a search by computer; the results are presented for links involving
up to 11 vertices (i.e. weight in the range 1 to 9). The same search
also yields the complete set of minimal strong links involving up to
10 vertices. We then analyse
the structure of the links by enquiring into their {\em reducibility}.
That is, we would like to know whether a given link can be understood as
a set of smaller links, or more generally as a set of smaller problems of
some kind. We consider two reductions which are already known, and we
present a third, called here the `short-cut'. The main result is
that the structure of links in general is of increasing complexity
as the size grows. However, analysis of games such as Hex can profit
from understanding the simplest classes, which are those which reduce
to one of the smallest irreducible links. These are all exhibited here.

\section{Definitions; notation}

The notation follows that of \cite{12steane}.
The Shannon game $(G,s,t)$ is played on a simple graph $G$
(no self-loops, no multiple edges)
with two vertices $s,t$ designated as {\em terminals}. We define the {\em size} of
a Shannon game to be the number of vertices in the graph (this is two
more than the number of vertices in the play area). A Shannon game
can also be considered to be a {\em link}. It is a strong link if \short~has
a second-player winnning strategy. It is a weak link if both players have
a first-player winning strategy.

A {\em pivot} of a weak link is a vertex which, if it is shorted, strengthens
the link (turning weak into strong).

A link of any kind is called {\em minimal} if and only if
the deletion of a single edge suffices to weaken the link.

An {\em induced Shannon game} is constructed as follows.
Let $(G,s,t)$ be a Shannon game and let $U$ be a set of vertices in $G$
having a neighbourhood $\Gamma(U)$ of size 2. Let $G'$ be the vertex-induced
sub-graph of $G$ induced by $U \cup \Gamma(U)$. Then $(G', \Gamma(U))$ is
an induced Shannon game having terminals $\Gamma(U)$ and play-area $U$
(it is a sub-game of $(G,t)$).
 
An {\em induced multi-Shannon game} is constructed
similarly, except that the neighbourhood of $U$ can have any 
number of vertices, resulting in any number of terminals in the induced
game.

A graph $G$ is said to {\em support} a Shannon game $(G,t)$ if the game can
be created on the graph by assigning two vertices in $G$ to act as terminals.

For the purposes of this paper,
a link will be said to {\em contain} a smaller link if the smaller link extends
between the same terminals as the larger, and the graph of the smaller link is a
sub-graph (not necessarily vertex-induced) of the larger link. For example,
S2 contains W1, but W3 does not contain W1.

We will make use of several results from 
`Threat, support and dead edges in the Shannon game'
\cite{12steane}. The notation `theorem TSD.n'
refers to theorem $n$ from that paper. We will also use the concepts of
{\em threat, supporting vertex, surrounding vertex, transverse edge,  
dead edge} from that paper.

\section{Minimal links}  \label{s.minlink}

We investigate the problem of links in general by first identifying
various generic properties, and then carrying out
an exhaustive search. 

\begin{table}
\begin{tabular}{r|rrrcc}
$n$ & connected & games   & sieved & minimal & S-irreducible\\
    &  graphs   &         & graphs & links   & links \\
\hline
 2 &          1 &       1 &      0 &    0    & 0 \\
 3 &          2 &       3 &      1 &    1    & 1 \\
 4 &          6 &      16 &      0 &    0    & 0 \\
 5 &         21 &      98 &      1 &    1    & 0 \\
 6 &        112 &     879 &      0 &    0    & 0 \\
 7 &        853 &   11260 &      9 &    5    & 1 \\
 8 &      11117 &  230505 &     35 &    0    & 0 \\
 9 &     261080 & 7949596 &    737 &   36    & 8 \\
10 &   11716571 & $\sim 5 \times 10^8$  &  21523 &   24    & 7  \\
11 & 1006700565 & $\sim 5 \times 10^{10}$ & O$(10^6)$ & 953 & 312
\end{tabular}
\caption{Statistics relating to the search for minimal weak links. The 
columns of the table are 1: size of graph; 2: number of connected graphs; 
3: number of non-isomorphic connected Shannon games; 4: number of graphs retained
by the sieve described in section \protect\ref{s.minlink}; 5: number of minimal weak
links; 6: number of Shannon-irreducible weak links. }
\label{t.numlink}
\end{table}

The number of non-isomorphic Shannon games on some connected graph
of size $n$ is given for small $n$ in column 3 of table \ref{t.numlink}.
This is the number of non-isomorphic ways of colouring two vertices in
some connected graph of size $n$ (we calculated it by generating and
colouring all such graphs, and checking them for isomorphism).
The number of these games that are weak links (i.e. those for which
\short\ has a 1st- but not a 2nd-player winning strategy) is unknown, but
it is certainly large. In order to get some general insight into the structure
of links, we focus our attention on {\em minimal weak links}. Obviously
any weak link contains one or more minimal weak links, and in order to find
the winner and a winning strategy for a game, it suffices to find only minimal
weak links. 

Once one has the set of
all minimal weak links of weight $w$, one can immediately deduce the
set of all minimal strong links of smaller weight, since they all appear
as {\em blocks} in those minimal links which have an articulation vertex. 
For example, a minimal weak link in which one terminal is a pendant
(a vertex of degree 1) yields a minimal strong link between the
neighbour of the pendant and the other terminal. Thus S2 can be
discovered by eliminating the pendant from W3, for example, and the two
non-isomorphic minimal strong links of weight 4 can be discovered by elliminating
the pendant from W5A and W5C. (The strong link obtained from W5B is isomorphic
to that obtained from W5A). 

In order to find these links, 
the concept behind our adopted method is to consider, in
principle, all connected graphs of
size $n$, and all possible terminal assignments for each graph,
but pruning the search as efficiently as possible. Before even
considering the terminal assignments, many graphs can be ruled out
because no terminal assignment in that graph will produce a minimal link.

\begin{lemma} \label{lem.dead}
Minimal weak links do not contain dead edges and do not contain (as a sub-graph)
any smaller weak link between the same terminals.
\end{lemma}
Proof. This is obvious; if there were a smaller link then the link in question
would not be minimal. \QED

\begin{lemma}  \label{lem.pendant}
A minimal weak link of weight $w>1$ can have at most one pendant vertex, and that vertex
must be a terminal.
\end{lemma}
Proof. 
A non-terminal pendant is dead so a link with such a vertex is not minimal. 
If both terminals have degree 1 then \cut\ has a 2nd player
winning strategy, namely to cut the neighbour of the remaining
pendant terminal after \short's first move, so we don't have a weak link.
This strategy must be available because the terminals are not both
adjacent to the same vertex, or the link would contain W1 and consequently
could not be minimal.  \QED

Next we will obtain some results which constrain the degree sequence of the
graph.  
Let the terminal vertices be $s$, $t$, and let their degrees be
$d_s$, $d_t$, and their neighbourhoods be 
$B_s \equiv \Gamma(s)$, $B_t \equiv \Gamma(t)$.
We will call these terminal neighbourhoods {\em borders}. Let the
vertices not in either border, and not terminal, be called the {\em centre}
of the link, and let $p$ be the number of vertices in the centre. Then, for
a link of weight $w$ in which the borders do not overlap (which is the case
for all minimal weak links except W1),
\be
w = d_s + d_t + p .    \label{sumdp}
\ee

\begin{theorem} \label{th.dterm}
For a minimal weak link of weight $w$, the size of the centre is constrained
as follows:
\be
\begin{array}{rcl}
w > 3 &\Rightarrow&  p \; \ge \; 1 \\ 
w > 5 &\Rightarrow&  p \; \ge \; 2
\end{array}
\ee
and for $w \le 7$ the terminal degrees $d_s$, $d_t$ satisfy
\be
d_s + d_t \le \frac{w+3}{2}.     \label{eq_dsdt}
\ee
\end{theorem}
Proof. The terminals are not adjacent in a weak link. Hence neither border
contains a terminal. Also, the borders do not overlap, for $w>1$,
or the link would contain W1 and therefore would not be minimal. 
To prove the first condition, suppose
the converse, i.e. $p=0$. Now, there can be no edges between a neighbour
of $s$ and more than one neighbour of $t$ or the link would contain W3
(see figure \ref{f.minweak}). But this implies that, after deleting
dead edges among each border, such
a link can only have the following structure:

\begin{picture}(40,35)(30,-5)
\multiput(47.5,7.5)(120,0){2}{\rule{5pt}{5pt}}
\put(52.5,11){\line(3,1){25}}
\put(52.5,9 ){\line(3,-1){25}}
\multiput(80,20)(0,-21){2}
{ 
  \multiput(0,0)(60,0){2}{\circle{5}}
  \put(2.5,0){\line(1,0){55}}
}
\put(142.5,-1){\line(3,1){25}}
\put(142.5,19){\line(3,-1){25}}
\end{picture}\\
with possibly more pairs of vertices similarly connected. Therefore the link
contains one or more mutually threatening pairs so is not minimal. (Also it is not a
weak link since \cut\ has a 2nd-player winning strategy.) We have a contradiction,
which completes the proof. The proof
of the second case ($w > 5$) uses a similar approach but is longer
and can be found in the appendix. The statement about $d_s + d_t$ follows
immediately by simple algebra using eqn (\ref{sumdp}) (and by using the
fact that we know it is
obeyed by W1 and W3, the only minimal weak links of weight $w \le 3$). \QED


\begin{theorem} \label{th.demax}
Every border vertex $b$ in a minimal weak link of weight $w$
has degree constrained by $d_b \le p + 2$, for $w>3$ and
$d_b \le p+1$ for $w>5$. ($p$ is the size of the centre
of the link, given by (\ref{sumdp})).
\end{theorem}
Proof. The vertex $b$ is by definition adjacent to one terminal. It cannot
be adjacent to the other terminal or the link would contain W1 so would not
be minimal. Also, $b$ cannot be adjacent to any other neighbour
of its terminal, or there would be an edge `across' the terminal, which would be a
dead edge (by theorem TSD.15), and obviously $b$ is not adjacent to itself.
Also, it can be adjacent to at most one
of the neighbours of the other terminal, otherwise the link would contain W3.
The result for $w>3$ follows. The proof for the case $w>5$
is given in the appendix.  \QED

\begin{theorem}  \label{th.dmax}
A minimal weak link of weight $w > 5$ satisfies
$d_{\rm max} + \min(d_s,d_t) \le w-1$, with one exception.
\end{theorem}
Proof. Consider an arbitrary vertex $v$. 
$v$ is either a terminal, or a border vertex, or a centre vertex.
If it is a terminal then its degree satisfies (\ref{sumdp}), 
hence $d_v + d_{\rm other} = w-p$ where $d_{\rm other}$ is the degree
of the other terminal. Now $\min(d_s,d_t) \le d_{\rm other}$ hence
we have, for this case,
\be
d_v + \min(d_s,d_t) \le w  - p .  \label{condterm}
\ee
If $v$ is a border vertex, then 
$d_v \le w - (d_s + d_t) + 1$ 
by theorem \ref{th.demax} and eqn (\ref{sumdp}). But
$d_s + d_t \ge \min(d_s,d_t) + 2$ since the terminals cannot both be
pendant. Hence the degree of a border vertex satisfies
\be
d_v + \min(d_s,d_t) \le w  - 1 .  \label{condbord}
\ee
If $v$ is neither a terminal nor a neighbour to a terminal (i.e. it
is a centre vertex), we show in the appendix that, with the
exception of the minimal link shown in figure \ref{f.exlink},
$d_v \le p-1 + \max(d_s,d_t)$. That is, the vertex cannot be adjacent
to all the largest border and all the rest of the centre and
a vertex in the other border.
Using eqn (\ref{sumdp}) we find condition (\ref{condbord}) again.
Noting from theorem \ref{th.dterm} that $p \ge 2$, we have that all
vertices in the graph satisfy (\ref{condbord}). In particular, the vertex
of maximum degree satisfies it. \QED

\begin{figure}
\myfig{0.15}{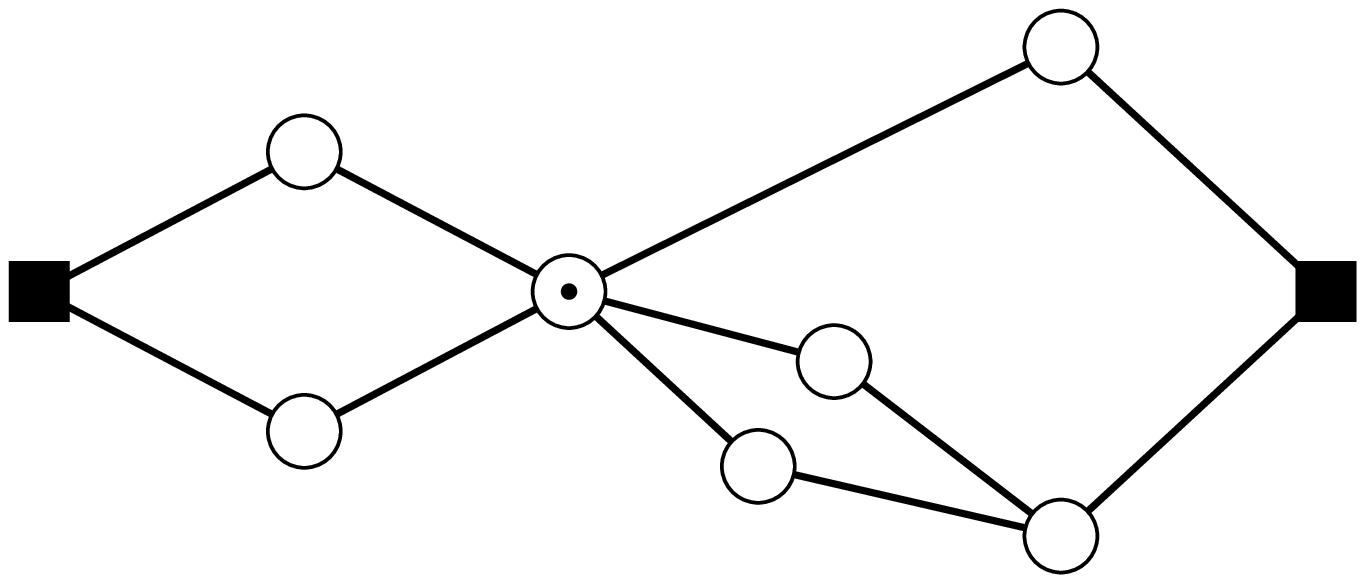}{f.exlink}
{The exception to theorem \protect\ref{th.dmax}.}
\end{figure}

Theorems \ref{th.dterm}, \ref{th.demax}, \ref{th.dmax} directly concern minimal links,
but they can also be used to constrain the set of graphs which can support a minimal link.
Consider an arbitrary graph on $n$ vertices, and let its degree sequence be ordered
such that $d_1 = d_{\rm min}$ is the smallest and $d_n$ the largest vertex degree. 
Then, for {\em any} choice of terminals $s,t$ on the graph we must 
have $d_1 + d_2 \le d_s + d_t$, and $d_1 + d_2 + d_3 \le d_s + d_t + d_b$.
Hence we have the following constraints from theorems \ref{th.dterm}, \ref{th.demax}, \ref{th.dmax}:
\be
d_1 + d_2 &\le \; \frac{n+1}{2} & \mbox{for}\; n \le 9,   \label{d1d2limA} \\
d_1 + d_2 &\le \; n - 4         & \mbox{for}\; n > 7    , \label{d1d2limB}  \\
d_1 + d_2 + d_3 &\le \; n       & \mbox{for}\; n > 5 ,    \label{d1d2d3limA} \\
d_1 + d_2 + d_3 &\le \; n-1     & \mbox{for}\; n > 7 ,    \label{d1d2d3limB} \\
d_1 + d_n &\le \; n-3           & \mbox{for}\; n > 7 ,    \label{d1dnlim} 
\ee
with one exception for the last condition.
Any graph which does not satisfy these constraints cannot support any
minimal weak link. In addition, lemma \ref{lem.pendant} furnishes the lower bound
\be
d_1 + d_2 \ge 3  \;\;\; \mbox{for} \; n > 3 .   \label{d1d2ge3}
\ee
Using $d_{\rm min} \le d_1,d_2,d_3$ in conditions (\ref{d1d2limA}) and (\ref{d1d2d3limB})
we find $d_{\rm min} \le (n+1)/4$ for $n \le 9$ 
and $d_{\rm min} \le (n-1)/3$ for $n > 7$.

\begin{theorem} \label{th.one_cut}
A minimal weak link has at most one articulation vertex (cut-vertex).
\end{theorem}
Proof. If the link is minimal then its graph is connected. But a connected
graph with at least one cut-vertex has at least two leaf blocks. Any terminal-free
leaf block is dead, but minimal links contain no dead vertices,
so the only possible arrangement is that there are two leaf blocks, containing
one terminal each. Then either the terminal blocks share a vertex, in which
case the rest of the graph is dead (and the rest of the graph is not empty under
the supposition that there is another cut-vertex), or they do not. If the
terminal blocks do not share a vertex then
there is a second-player winning strategy for \cut, namely, to disconnect (by deleting
the appropriate cut-vertex) whichever
terminal block is not `rescued' by \short's first move, hence in this case 
there can be no weak link. \QED 

\begin{theorem} \label{th.bridge}
A link is a minimal weak link if and only if its bridge-reduced form is a
minimal weak link.
\end{theorem}
Proof.
The theorem concerns links which induce one or more `bridges', i.e.
the induced Shannon game S2. By the `bridge-reduced form' we mean the link
obtained after replacing the bridge by an edge---that is, if $a$,$b$ both
have degree two and the same neighbours, then we add an edge between the neighbours
and then delete $a$ and $b$.  If the added edge
is a spectator of the reduced game, then so was the bridge in
the original game. If not, then at some stage \cut\ has to try to prevent \short\
from acquiring the edge; the bridge is a minimal strong link that ensures \cut\
will not be able to do so, and it does not otherwise affect the game. 
(Note that, in view of theorem TDS.16, neither of the vertices in the bridge
is a terminal if the link is minimal. Also, if the pair of neighbours was already
adjacent then the added edge results in a non-simple graph (one having 
two edges between given vertices) which is clearly non-minimal.)
\QED

\begin{lemma} \label{lem.support}
In a graph supporting a minimal weak link other than W1, no vertex is in
more than one mutually surrounding pair.
\end{lemma}
Proof.
If a mutually surrounding pair are both terminal then we have either W1 or
a strong link such as S2. If one and only one of a 
mutually surrounding pair is terminal, then the other is dead by theorem TSD.16.
It follows that, in minimal weak links other than W1,
mutually surrounding pairs must be terminal-free.
Suppose $a,b$ are mutually surrounding and so are $a,c$. Then $c$ is redundant, that is,
$c$ and its edges can be deleted without changing the outcome of the game. This is
because a mutually surrounding pair such as $a,b$ is `captured' \cite{12steane}, meaning
both can be filled in by \short\ as a free move, without changing the outcome of the game.
After that filling in, $c$ becomes simplicial (its neighbours form a clique)
so it is dead (see \cite{06Hayward} and theorem TSD.13). \QED

\begin{lemma} \label{lem.trans}
In a minimal weak link other than W1 there can be at most one transverse edge,
and if there is such an edge then it is to a pendant terminal.
\end{lemma}
Proof.  We show that if the condition does not hold then there is a dead
edge, which is ruled out by lemma \ref{lem.dead}.
A transverse edge to a non-terminal vertex is dead. Hence if there is
a transverse edge in a minimal link, then it is to a terminal. Let the vertices
at the ends of this edge be $t,v$, with $t$ the terminal. By definition of
a transverse edge, either $t$ surrounds $v$ or $v$ surrounds $t$. But if
$t$ surrounds $v$ then $v$ is dead, by theorem TSD.16
(unless $v$ is also a terminal, but then we have adjacent terminals so not a weak link). 
If $v$ surrounds $t$ and $t$ has degree greater than one, then a neighbour of $t$
(namely, $v$) is adjacent to at least one other neighbour of $t$, so there is
a dead edge by theorem TSD.15. \QED

\begin{lemma} \label{lem.pivot}
Pivots of a minimal weak link are triangle-free, not threatened, not supported, and
do not support any non-terminal vertex.
\end{lemma}
Proof. The unthreatened property is obvious. Now
consider triangles. A pivot of a weak link is a vertex that can serve as the
opening move for \short\, such that he wins the game. Such a 
move by \short\ will in any case introduce edges between all the
neighbours of the pivot, so there was no need for those edges to be present
before the opening move. In the case of multiple pivots, the argument applies
to each one. Next consider support. If a pivot is supported, then the supporting
vertex must also be a pivot (the latter cannot be a terminal because then we would have
a dead vertex by theorem TDS.16). However, if a pivot supports another vertex, then
after an opening move at that pivot, the supported vertex will be dead (except
when it is a terminal), which implies it was not needed in any case and the link 
would not be minimal.  \QED

\begin{lemma} \label{shortcut_dead}
If a pivot $b$ of a minimal link is in border $B_s$, then
there are no edges between $\Gamma(b) \setminus \{ s \}$ and
$B_s \setminus \{ b \}$.
\end{lemma}

Proof. Let $c$ be a vertex in $\Gamma(b) \setminus \{s\}$, i.e.
it is a neighbour of $b$ other than the terminal. We know
that shorting $b$ is a possible opening move.
After this there is an edge $sc$. It follows that if there
was an edge between $c$ and  $B_s \setminus \{ b \}$ in the opening
position, then that edge is now dead by theorem TDS.15. Therefore
it must have been a spectator in the opening position (since we
can now delete it with impunity, and the opening move at $b$ is
certainly available as one strategy for \short\ to win). But
a minimal link has no spectator edges. \QED

\subsection{Sieving for minimal links}

The following list gives conditions that a {\em graph} must satisfy in order
for it to be possible that that graph might support a minimal link. 
That is, if there exists some terminal assignment giving a minimal link,
then the graph must satisfy the conditions. They are necessary but not sufficient
conditions for the existence of such a terminal assignment. In order to
search for minimal links, an algorithm that could in principle generate all
graphs of given size was used, but these conditions were implemented,
in the order given here, i.e. simplest first,
in order to rule out graphs or sets of graphs as quickly as possible.
This is a type of `sieve' which retains only graphs which might possibly
support a minimal weak link. 

\begin{enumerate}
\item Constraints (\ref{d1d2limA})--(\ref{d1dnlim}), (\ref{d1d2ge3}).
\item The graph is connected.
\item There is at most one transverse edge (lemma \ref{lem.trans}).
\item If there is a transverse edge, it is incident on a pendant (lemma \ref{lem.trans}).
\item There are at least 3 triangle-free vertices. Proof:
we require there to be two vertices that could serve as terminals, and a third which
could serve as a pivot. But terminal vertices are triangle-free in a minimal link,
because otherwise there would be a dead edge across a terminal (theorem TSD.15),
and pivots are triangle-free by lemma \ref{lem.pivot}.
\item The bridge-reduced graph might possibly support a minimal link
(theorem \ref{th.bridge}).
\item \label{cond.pendant}
If there is a pendant, then its neighbour is triangle-free and of degree $d>2$
when $w>1$. Proof:
if there is a pendant, it must be a terminal (lemma \ref{lem.trans}) and the
vertex $v$ adjacent to that terminal must be the only pivot, since otherwise,
after \short's first
move elsewhere, \cut\ can win the game by cutting $v$ and thus isolating the terminal.
Hence $v$ is triangle-free. It cannot be of degree 2 since then \cut\ can win by
cutting either it or its neighbour, which \short\ cannot prevent.
\item No vertex is in more than one mutually surrounding pair (lemma \ref{lem.support}).
\item Mutually surrounding vertices must have no triangles. Proof:
a mutually surrounding pair can be filled in at the outset without changing the game
outcome. Therefore any edges amongst their neighbours will be supplied in any case
and are not needed in the opening position.
\item \label{cond.terminals}
There is at least one candidate terminal pair.
Terminals of a minimal link are triangle-free (see above) and do not surround any
other vertex (theorem TSD.16). 
If there is a pendant, it must be one of the terminals.
Candidate terminal pairs are not adjacent and for $n>3$ have no 2-walk between them
(or the link would contain W1). The sum of their degrees must satisfy
theorem \ref{th.dterm}.
\item \label{cond.threat} No pair of non-terminal vertices may threaten each other.
Having identified candidate terminals, one can look among the other vertices
in the sure knowledge that none will be assigned as a terminal. If any of those
are mutually threatening then they are `lost' \cite{12steane}, that is, they can
be deleted without changing the outcome of the game. This would imply the link was
not minimal.
\item There is a candidate opening move.
Candidate opening moves must satisfy the conditions on a pivot given by
lemma \ref{lem.pivot}. Note that a threatened vertex
either has a triangle or is of degree $d=2$, so it is sufficient to check
for triangles and the degree here.
(This check can be omitted if there is a pendant, since
condition \ref{cond.pendant} already implies it.)
Also, if a vertex is known to be terminal (for example because it appears in all
candidate terminal pairs) then it cannot be a pivot of any weak link supported
by the graph. Candidate pivots must not
be surrounded, nor may they surround any vertex $v$ known to be 
non-terminal (this may be known, for example, because $v$ does not appear
in any candidate terminal pair).
\item If $n>5$ then for each candidate terminal pair, no neighbour of
one terminal should
have more than one 2-walk to the other terminal, or the link would contain W3.
Also, for $n>7$ no vertex should have more than one 2-walk to both terminals, or
the link would contain either W5X or a smaller link.
Also, the constraint described in theorem \ref{th.demax} must be satisfied.
\item For each candidate terminal pair, no other vertices should be in
a mutually threatening pair. This is like condition \ref{cond.threat}, but
having applied it in step \ref{cond.threat} to the set of vertices which 
are not terminal in any candidate link supported by the graph, we here apply
it again, on a case by case basis, for each available choice of a terminal pair.
This provides a modest tightening of the sieve at little calculational cost.
\end{enumerate}

\begin{table}
\begin{tabular}{r|rccccrrrrr}
$w$& links&$p_{\rm min}$& $E$ &$d_{\rm max}$&  S-    &  P-    &   T-   & SPT-  \\
   &      &             &     &             & irred. & irred. & irred. & irred. \\
\hline
 1 &    1 &          & 2--2   &  2--2      & 0     & 1     &  1 & 0 \\
 3 &    1 &    0     & 5--5   &  3--3      & 0     & 0     &  0 & 0 \\
 5 &    5 &    1     & 8--9   &  3--4      & 1     & 1     &  2 & 0 \\
 7 &   36 &    2     & 11--13 &  3--5      & 8     &   8   &  20 & 2 \\
 8 &   24 &    3     & 14--16 &  4--5      & 7     &  24   &  24 & 7 \\
 9 &  953 &    3     & 14--21 &  3--6      & 312   & 544   & 766 & 208
\end{tabular}
\caption{Some statistics of minimal weak links. The columns are
$w$: weight; {\em links}: number of minimal links; $p_{\rm min}$: smallest centre size; $E$: least and greatest
number of edges; $d_{\rm max}$: least and greatest maximal degree; last four columns: number of irreducible links.}
\label{t.weaklink}
\end{table}

\begin{table}
\begin{tabular}{r|rccccrrrrr}
$w$& links&$p_{\rm min}$& $E$ &$d_{\rm max}$&  S-    &  P-    &   T-   & SPT-  \\
   &      &             &     &             & irred. & irred. & irred. & irred. \\
\hline
 0 &    1  &          & 1--1   & 1--1    & 1     & 1     &  1 & 1 \\
 2 &    1  &          & 4--4   & 2--2    & 0     & 0     &  1 & 0 \\
 4 &    2  &   0      & 7--8   & 3--3    & 1     & 0     &  1 & 0 \\
 6 &   14  &   0      & 10--12 & 3--4    & 4     & 2     &  6 & 0 \\
 7 &   10  &   1      & 13--15 & 4--5    & 10    & 10    &  5 & 5 \\
 8 &  304  &   1      & 13--20 & 3--6    & 196   & 163   & 204 & 120
\end{tabular}
\caption{Statistics of minimal strong links; the column headings are the same as in 
table \protect\ref{t.weaklink}. }
\label{t.stronglink}
\end{table}

For $n <11$ we used the above conditions to sieve a list of
all graphs of given size. 
For $n=11$ we proceeded as follows.
First all graphs of size $n=9$ were
generated, and then a vertex added in order to form of list of graphs
on ten vertices, including disconnected ones, in which at least two
vertices are triangle-free and there is no pair of adjacent
mutually surrounding vertices. Next, a further vertex was added
in such a way as to respect all the above conditions, thus
generating a list of all such graphs on 11 vertices (including
some isomorphic duplicates).

The number of graphs retained is shown in table \ref{t.numlink}. 

One can also suggest further conditions to tighten the sieve.
For example, for a given terminal pair, one could use
lemma \ref{shortcut_dead} to further constrain the set of
candidate opening moves, or for $n > 7$ one could implement a check for the 
presence of any link of weight 5. Another idea
is to investigate the possible existence of a triangle-free vertex
when \short\ makes his second move, making due allowance for the
effects of the first two moves by \short\ and \cut. We investigated
a test for this, but it was found to offer no improvement in the sieve
for $n \le 9$.

For each retained graph, we next assigned
terminals in each distinct (i.e. non-isomorphic) way that might possibly give
a minimal weak link (i.e. the terminals are triangle-free, non adjacent, etc., as
in item \ref{cond.terminals} above), and for each case tested the resulting link.
We calculated the strength of the link,
using a modest Shannon-game-solving algorithm. Whenever a weak link was found
(i.e. one for which \short\ had a 1st- but not a 2nd-player winning strategy),
we tested whether it was minimal by deleting each edge in turn and re-solving
the game. We thus found all minimal weak links of given size. (A
strategy which speeds up the calculation for larger links is to first
check for the presence of an induced Shannon game, and check each such
game for minimality; another useful strategy is to test the deletion of
triangle edges first, since these are more likely to be spectators.). 

The number of minimal weak links for small values of $n$, calculated by this
method, is shown in table \ref{t.numlink}, and table \ref{t.weaklink} gives
some further information which will be discussed in the following sections. 
Figure \ref{f.minweak} shows all the minimal weak links
of weight up to 5 (i.e. $n \le 7$); further figures will be discussed below.

The set of minimal strong links of weight $w$ was obtained from the set of 
weak links of weight $w+1$ by deleting the pendant from those that have
a pendant (as discussed at the beginning of this section) and performing
an isomorphism check, so as to retain only one example from each isomorphism
class. The numbers found, and some other properties, 
are given in table \ref{t.stronglink}.

For small $n$, it is easy to satisfy oneself that any given
link found by this method is weak, or strong, as the case may be, and
minimal. It is, however, hard to prove that the result is correct, that is,
that there exist no other minimal links, because the method relies on
a computer calculation to handle the large number of candidate graphs and
terminal assignments. As a modest step towards a proof of correctness,
the appendix presents a proof that there is no minimal weak link of weight 4.

The absence of minimal weak links of weights 2,4,6 suggests the question, are
there are any minimal weak links of even weight? This is answered in the affirmative
by our search, which finds 24 such links of weight 8. It immediately follows (by
replacing edges by bridges) that
there are minimal weak links of all higher weights, whether even or odd.
We also find that there
exist minimal weak links having one or more triangles; the smallest of these
are of weight 8. We will discuss these results further after first
investigating the concept of reducibility.

\section{Reducibility}

For any given weak link one can always obtain a larger weak link by replacing
one of its edges by a strong link such as S2. For example, W1 can be expanded to W3,
and, conversely, W3 can be {\em reduced} to W1 by replacing the bridge by
an edge.
In order to explore weak links in general, it suffices to list, or otherwise
understand, just those that can not be reduced in this way or other simple ways. 

We will consider only minimal weak links. Obviously a non-minimal link can be
`reduced' in the sense of deleting `spectator' edges (those whose deletion,
in a given link, does not weaken the link), and often this can be done in more than
one way. However, identifying such edges
is itself a non-trivial task. We will restrict attention to minimal links
not because that is necessarily the best way to solve the Shannon game, but
because it is one natural way to get insight into the general
problem, and is an interesting area in its own right.

A minimal link can be said to be `reducible' if there is some way of
dividing it up into sub-problems. The example of W3 noted above is one
simple reduction, but a general notion of reducibility is not easy to
formulate. One idea that suggests itself is to call a link reducible if
it induces a smaller multi-Shannon game with a definite outcome,
i.e. one which is won outright by one player or the other, not drawn.
It could, for example, be captured by \short\ or dominated by both players
(i.e. both have a first-player winning strategy). In the first case
one can fill in the play area of the induced game with \short's
colour without changing the outcome.
In the second case one can replace the
play area of the induced game by a single vertex suitably connected to
the terminals. However, it is hard to identify such games, and this is
not the only type of reduction of interest
(we will consider another type below). To discover
whether a multi-Shannon game is not drawn is a smaller task than both to 
discover that and also identify the winner, but it remains a hard
task in general. However, an important exception is the case where there are just two
terminals---i.e. when the induced multi-Shannon game is in fact a Shannon game. In
this case one knows immediately that there is a winner, since Shannon games
cannot be drawn. Our strategy in the following will be to exploit this
and two other reductions
which are easy to identify in any given link. This does not rule out
the possible existence of further easily-identifiable reductions.


{\bf Reductions.} Consider the following properties that a link may possess.
\begin{enumerate}
\item
It induces a smaller Shannon game of weight $w>1$.
\item  
It induces a multi-Shannon game with a two-vertex play area having a 2nd-player
winning strategy for \short.
\item
It induces a {\em short-cut}. A {\em short-cut} is a pair of adjacent
non-terminal vertices, one of which has degree 2 and is adjacent to one
terminal and at distance 2 from the other terminal.
\end{enumerate}
A link having these properties will be said to be {\em Shannon-reducible},
{\em pair-reducible} or {\em short-cut-reducible}, respectively, and these
will be abbreviated to S-, P-, T-reducible. W3 is `SPT-reducible' because it
has all three properties. W5X is SP-reducible because it has properties 1 and 2,
but it is not T-reducible. W5Z is T-reducible but not S- or P-reducible.

The above reductions are useful because their presence is easily detected.
To discover whether a link induces a non-trivial Shannon game, use the
fact that the play area of an induced Shannon game has a cutset
of size 2. So it suffices to identify sets of vertices having a cutset of size 2.
This can be done, for example, by deleting each vertex in turn and finding
articulation vertices (cut-vertices) of the resulting graph using a depth-first-search.

Case (2) (multi-Shannon game on two vertices won by \short) can be detected 
using the concept of {\em support} as discussed in \cite{12steane}. 
Let the play area vertices be $a,b$. If \short\ has a 2nd-player
winning strategy then, by definition of `win' for the multi-Shannon game,
$G \ast a - b = G - a \ast b = G \ast a \ast b$
(in the notation of \cite{12steane}). Therefore $a$ and $b$ `support' one another.
To find such mutually supporting pairs,
one may use the concepts of {\em dead edges} and {\em surrounding} vertices
described in \cite{12steane}. First, one identifies every vertex that surrounds another,
for example by using the fact that if $A^2(u,v) = d_v - A(u,v)$ then $u$ surrounds $v$,
where $A$ is the adjacency matrix and $d$ is the vertex degree. Next, one deletes
`transverse' edges (those between a vertex and a surrounding vertex) where neither vertex
is terminal, because such edges are dead. After doing this, mutually supporting pairs
are also mutually surrounding pairs, which can be found by checking for pairs that
satisfy $A^2(u,v) = d_v = d_u$.

Case (3) is easily detected, either by examining short walks from all 
vertices of degree 2, or
by examining the neighbourhoods of the terminals. However, we need to
prove that this case leads to a reduction. A short-cut is neither a Shannon game,
nor a multi-Shannon game won outright by either player. Instead, it is a simple
combination of vertices in which one can be shorted, the other cut, without changing
the outcome of the larger game, as we shall now prove.

\begin{theorem}
In any minimal link, or non-minimal weak link with \short\ to play,
and in any (not necessarily minimal) Shannon game 
in which either player has a 2nd-player winning strategy, 
the outcome of the game is not changed if the vertex of larger degree 
in any given short-cut is shorted, and the other is cut.  
\end{theorem}
Corollary. If a weak link has a short-cut, then the $d>2$ vertex in the short-cut
is a pivot of the link.

Proof. Let $\cal S$ be the Shannon game in question, and let
the vertices in the short-cut be $a,b$, with $d_b = 2$ and
$d_a \ge 2$ (the definition of a short-cut implies that one vertex has 
degree 2 and the other has degree $\ge 2$).
If $d_a=2$ then $a,b$ are mutually threatening and therefore lost,
see \cite{12steane}. This means they can both be cut without changing the outcome
of the larger game, but this is equivalent to shorting one and cutting the other
since for a mutually threatening pair, $G \ast a - b = G - a \ast b = G-a-b$.
This completes the proof for the case $d_a=2$
(this case does not arise in a minimal link). 

If $d_a > 2$, we consider the following cases:
(A) \cut\ has a 2nd-player winning strategy for $\cal S$,
(B) $\cal S$ is a weak link,
(C) $\cal S$ is a strong link.
In case (A) we can allow
\short\ to move first without changing the outcome, and all
opening moves by \short\ are losing. \short\ could, for example, 
open at $a$. \cut\ must then reply at $b$, or he will lose. This position must
still be won by \cut\ since \short\ can force it from the opening position, which
was claimed to be won by \cut. This completes the proof for case (A).

In case (B), first suppose that \short\ is to move first. 
\short's winning strategy must require him to take possession
of either both $a$ and $b$, or one of them, or neither. It cannot require
both, because \cut\ can certainly prevent that by cutting one of them
on \cut's first move. It cannot require $b$ alone, because $b$ is
threatened---if \short\ does not open at $b$ then \cut\ can cut it; if \short\
does open at $b$ then a reply at $a$ produces a position 
equivalent to having cut both of them. Therefore \short\ requires either $a$
but not $b$, or neither of them. If he requires $a$ then he can safely open there,
because \cut's reply is forced---it must be $b$ or he will lose immediately.
\short\ then has a position in which two vertices have been occupied as required
in his winning strategy, or else in which he has gained a vertex ($a$) which he
did not require, without the loss of any vertex which he did require. Owning
an extra vertex is never detrimental in the Shannon game, so in either
case the position remains a won position for \short. This completes the proof
of both the theorem and the corollary in the case where \short\ opens. 
Next, still in case (B), suppose that \cut\ opens.
The theorem only concerns minimal links in this situation. 
We need to prove that \cut\ wins the reduced game. Let's assume the contrary.
Then, after the reduction, all moves by \cut\ are losing, i.e.
$G \ast a - b - v$ is a weak link for all $v \ne a,b$. But this implies that
$(G-v)$ must itself be a weak link. For if \short\ were to open on $(G-v)$ 
at $a$ he would force a reply at $b$ and thus
attain the position $G - v \ast a - b = G \ast a - b - v$ which he wins.
Hence $(G-v)$ is a weak link, which implies $G$ is not a minimal
weak link, which contradicts the premise. 

Finally, consider case (C). We will prove that $G \ast a - b$ is a strong link.
If $G$ is a strong link then all opening moves on $G$ 
by \cut\ are losing; in particular, $G-v$ must be a weak link, where $v \ne a,b$. Now, 
$\{a,b\}$ remains a short-cut on $G-v$, so we can apply the result already obtained
for \short\ openings on a weak link, and deduce that $G-v \ast a-b$, with \short\
to play, is won by \short. Hence $(G \ast a - b)-v$, with \short\ to play, is
won by \short\ for all $v \ne a,b$. Hence $(G \ast a - b)$ is a strong link. \QED

\section{Discussion}

We first discuss some structures of interest, and then some general trends.

W5Z is interesting for several reasons. It is the smallest minimal link to
have more than one pivot, and it is the smallest Shannon-irreducible weak
link (and it is also P-irreducible).
This means it will never be discovered by an algorithm such
as `H-search' \cite{00Anshelevich}
which uses only the AND, OR rules. However, it is sufficiently small to be reasonably
common and worth the effort of detecting. In the game of Hex, for example,
its detection, combined with the OR-rule and various S-reducible links, 
suffices to allow a number of useful `edge-templates' to be constructed.

We noted in table \ref{t.numlink} that there are 36 minimal weak links of
weight 7. Of these 28 are S-reducible; these may be regarded
as the more `natural' enlargements of one or more of the minimal links of weight 5. 
The 8 S-irreducible weak links of weight 7 are shown in figure
\ref{f.n9}. 2 of these are P-reducible and not T-reducible, 3 are 
T-reducible and not P-reducible, 1 is PT-reducible, and 2 are
SPT-irreducible. The P- or T-reducible links are reasonably obvious enlargements
of one or more of the W5 links. 
The diagrams of W7A, W7B, W7C in figure \ref{f.n9} have been drawn 
in such a way as to make it easy
to see which W5 link is obtained by the P-reduction, i.e. by shorting the 
pair of mutually supporting
vertices. W7A reduces to W5X, W7B and W7C reduce to W5Z (with an extra edge
in each case). Working in the opposite direction, for any weight $w$ one can
obtain a P-reducible W$w$ link by replacing two or more edges in a W$(w-2)$ link
by a mutually supporting pair of vertices suitably joined to the
vertices at the ends of the replaced edges. One can obtain a T-reducible W$w$ link 
from a W$(w-2)$ link by inserting a vertex $a$ into an edge from one terminal, 
and adding another new vertex $b$, adjacent to $a$ and the other terminal,
such that $a,b$ form a short-cut with $b$ of degree 2.

There is one example at weight 7 of a graph which supports more than one
minimal weak link. This may be understood as a result of the fact that
in the graph of W5X, there are two different pairs of vertices that can
serve as terminals and give a minimal link. These two terminal
assignments are isomorphic in W5X, but when an edge between a $d=2$
vertex and a terminal in W5X is replaced
by a bridge, the two terminal assignments in the resulting weight 7 link
are not isomorphic.

\begin{figure}
\myfig{0.6}{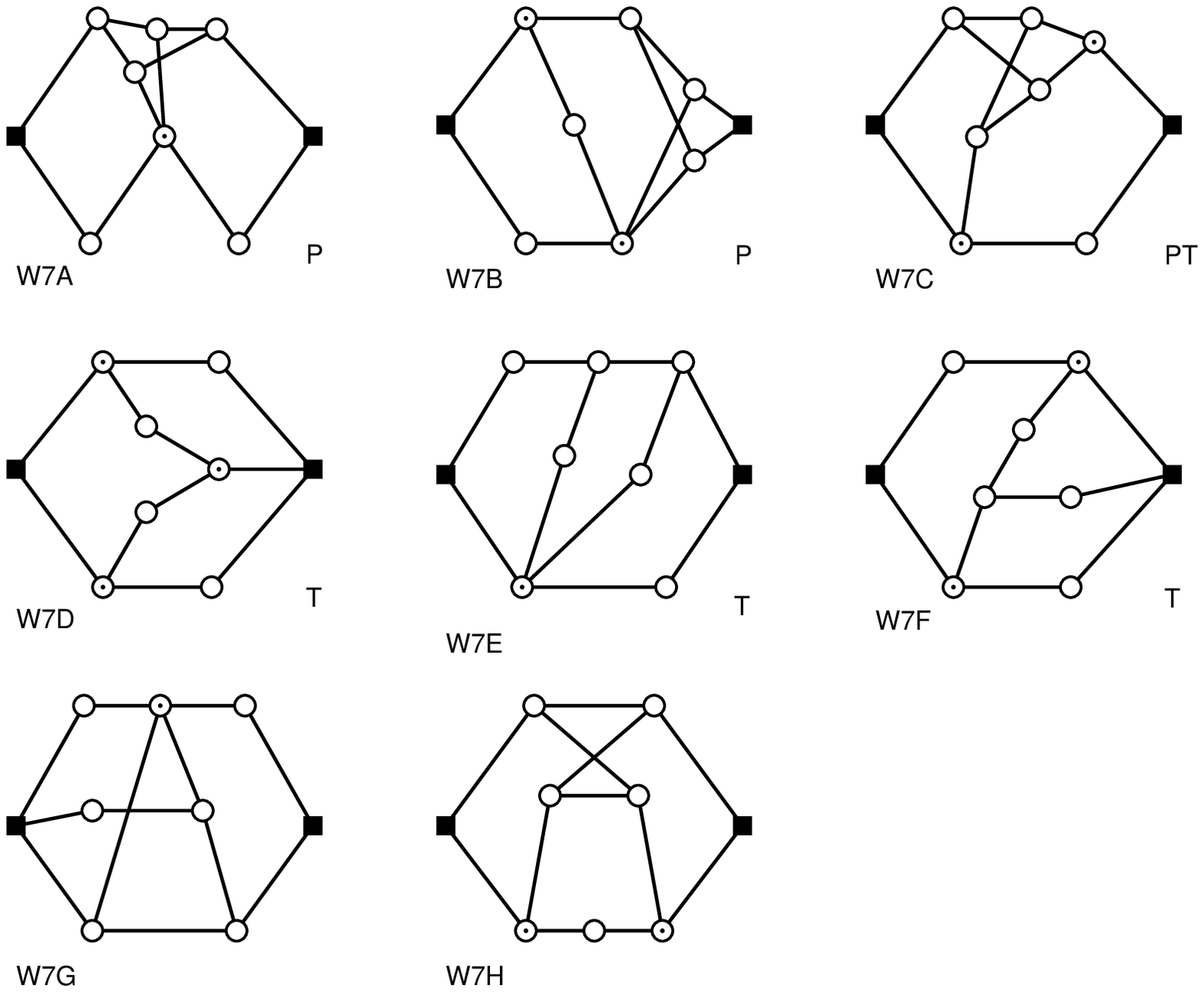}{f.n9}
{All Shannon-irreducible minimal weak links of weight 7.
The labels `P', `C', `PC' indicate those which
are P-reducible, T-reducible and PT-reducible, respectively.
Dotted vertices are pivots.}
\end{figure}

\begin{figure}
\myfig{0.18}{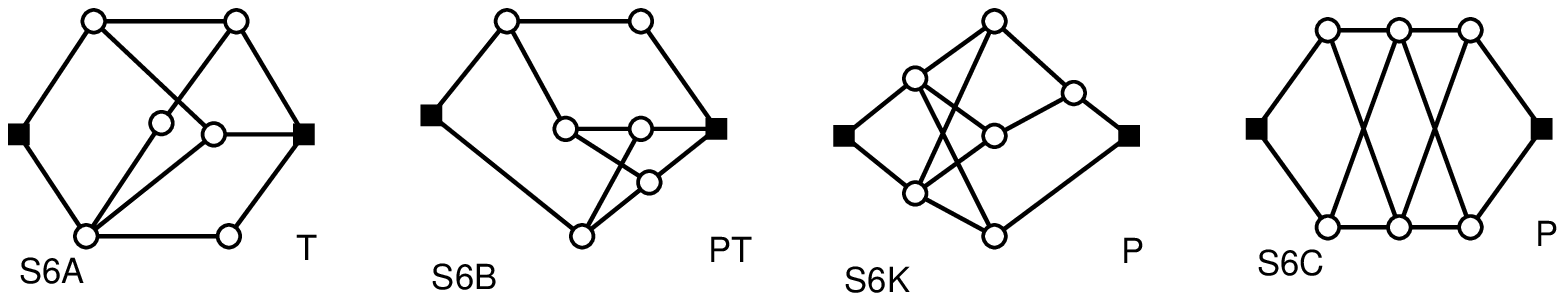}{f.s8}
{All Shannon-irreducible minimal strong links of weight 6.
The set of overlapping weak links forming each strong link is as follows,
where a superscript index indicates the number of each type of weak link required:
S6A: \{W3$^2$, W5Z\}; S6B: \{W3, W5B, W5Z$^2$\};
S6K: \{W3, W5C, W5A$^2$\}; S6C: \{W5C$^4$, W5X$^2$\}.
}
\end{figure}

W7F is an example of a situation known in Hex is a `forking ladder escape'.
The lower pivot is at the foot of a short two-rung `ladder'. In the absence
of the upper pivot, \cut\ could prevent \short\ from making this ladder connect
to the right terminal. The upper pivot is a `ladder escape', and furthermore
it is a `forking ladder escape' because by occupying it, \short\ both
gains the ladder escape he wants, and also forces \cut\ to deal with the
other winning link he is threatening to complete. W7E is an example of
another type of `ladder', one which runs longitudinally (i.e. from one
terminal to the other) rather than transversely (i.e. across the neighbourhood
of one terminal) as W7F.

The two SPT-irreducible W7 links represent `new' structures which cannot be 
obtained in a simple way from any W5 link. They do not
have any readily calculable properties that would make them easy to detect
in a larger Shannon game.

There are 24 minimal weak links of weight 8.
None of these are P-reducible, and none are T-reducible. 
This is related to the fact that there are no minimal weak links of weight 6.
There could only be a P- or T-reduction of a W8 link if the resulting link
(of weight 6) had a spectator vertex. However, in order for that to occur it
seems that the original W8 link must have been non-minimal. This suggests
the following, which I leave as a conjecture:

{\bf Conjecture}. If a minimal weak link of weight $w$ is either P-reducible
or T-reducible, then there exists a minimal weak link of weight $w-2$.

17 of the W8 links are S-reducible. Of these, 3 have no triangles and 14 have two triangles
(none have just one triangle). 7 of the W8 links are S-irreducible; these are 
shown in figure \ref{f.n10}. Of these, 3 have no triangles and 4 have one triangle.

At weight 9 there are too many minimal links to show them all; figure \ref{f.n11}
shows a few specimens having interesting properties.

\begin{figure}
\myfig{0.4}{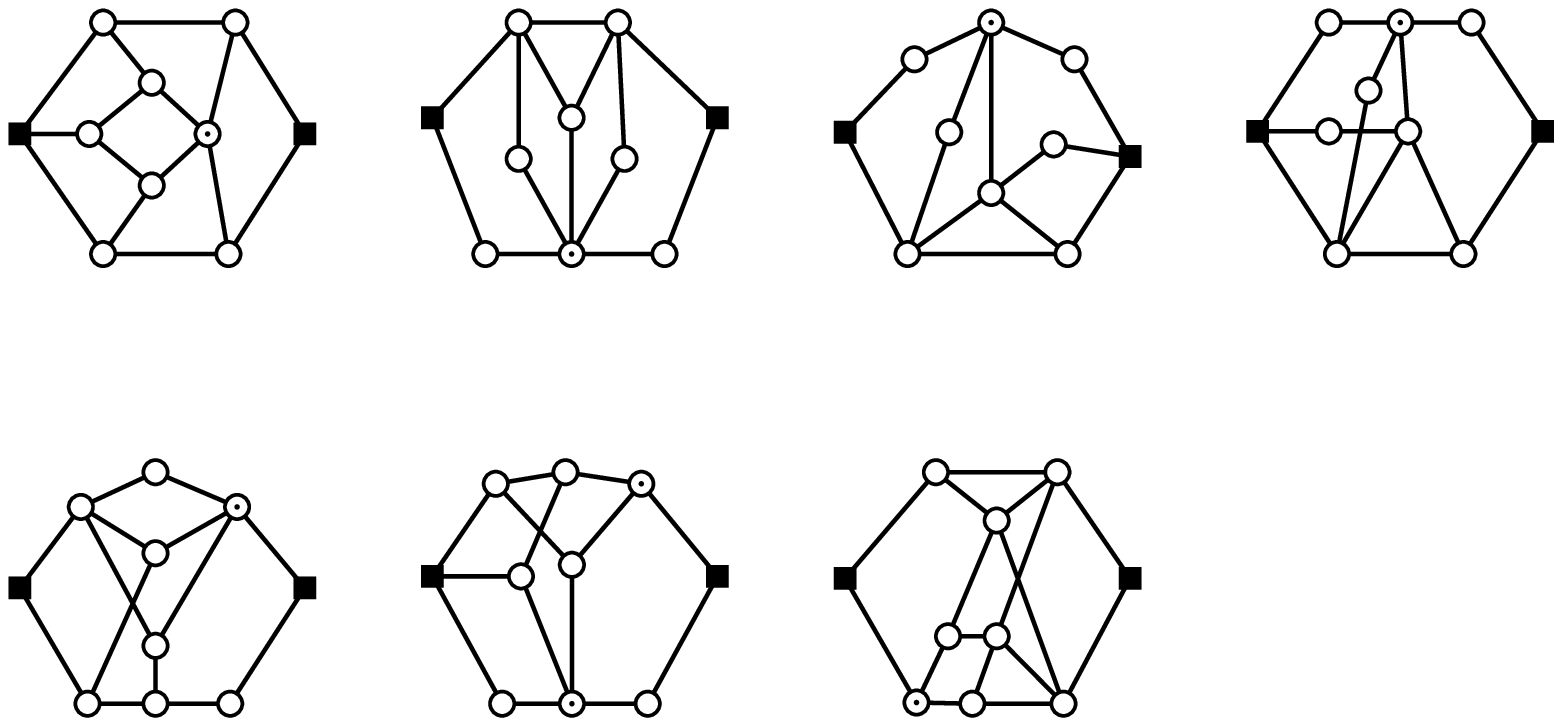}{f.n10}
{All Shannon-irreducible minimal weak links of weight 8.
Dotted vertices are pivots.}
\end{figure}

\begin{figure}
\myfig{0.4}{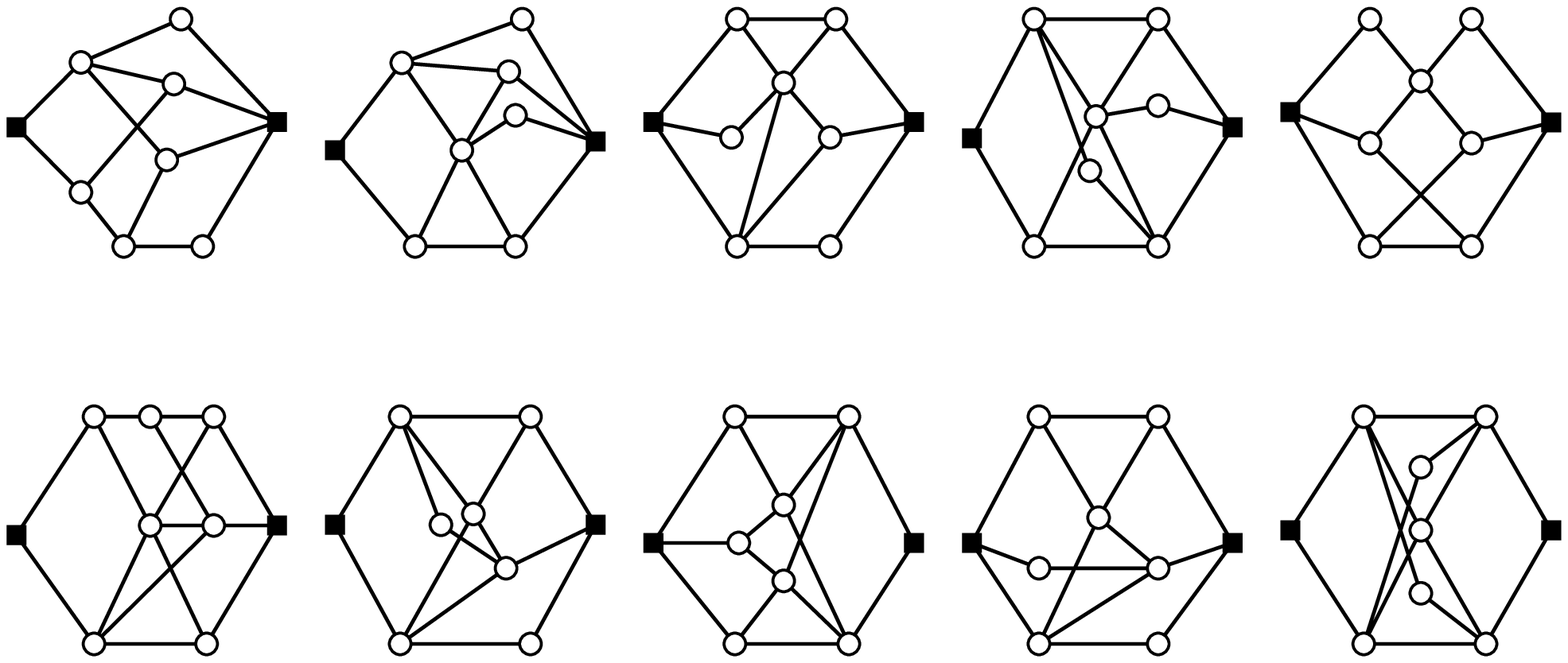}{f.n9s}
{All Shannon-irreducible minimal strong links of weight 7.
The decomposition into overlapping weak links is not always unique for any given
strong link; the following lists one possibility for each of the above:
\{W3$^2$, W5B\}; \{W3, W5B, W5X\}; \{W3, W5x, W5Z\};
\{W5X$^3$, W5B$^2$, W5Z\}; \{W3$^2$, W5X\}; \{W3, W5B, W5C, W5X\};
\{W3, W5A, W5Z\}; \{W5A$^2$, W5Z$^2$\}; \{W3$^2$, W5Z$^2$\}; \{W5A$^4$, W5Z\}. 
}
\end{figure}

\begin{figure}
\myfig{0.25}{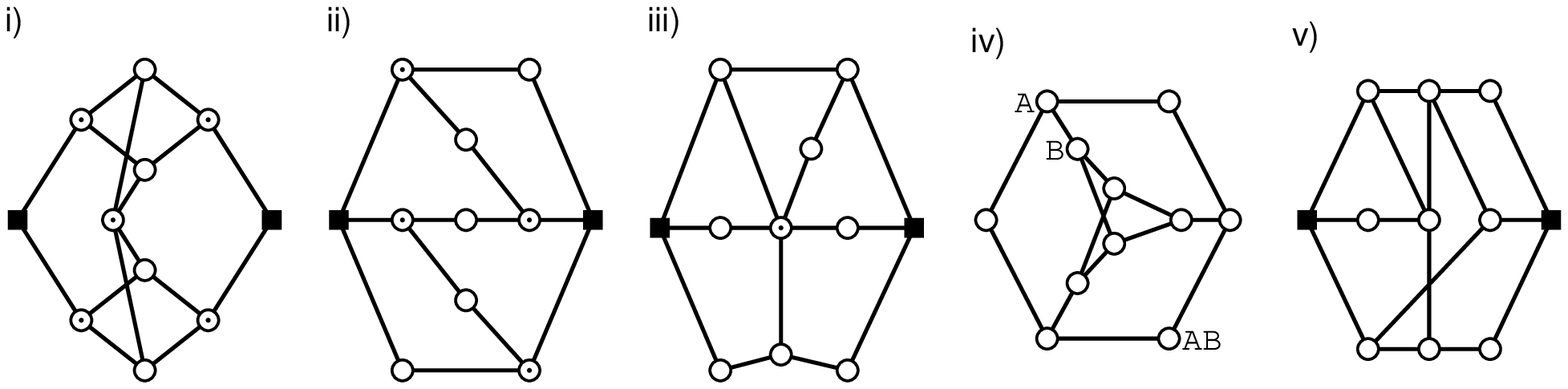}{f.n11}
{A few interesting minimal weak links of weight 9.
(i) and (ii) are the only ones having more than 3 pivots; (iii) is the only S-irreducible
one that meets the bound of theorem \protect\ref{th.dmax}; (iv) has the largest
set of different terminal assignments that lead to minimal weak links---the terminal
AB is minimally weakly linked to either A or B,
and owing to the 3-fold automorphism (which can be seen as a 3-fold rotational symmetry
of an embedding of the graph in three-dimensional space) there
are a total of 6 possible terminal pairs which are minimally weakly linked;
(v) is an example link having no simple reduction; the reader is invited to
find the pivot.}
\end{figure}

\subsection{Limits on degree and distance}

Table \ref{t.weaklink} allows a few general properties to be discerned.
The number of minimal weak links grows rapidly with $w$; the limited data
are consistent with a growth roughly in proportion to the square root of
the number of graphs on $w$ vertices. 
The size $p$ of the centre of the link is consistent with theorem 1 (of course),
and the fact that we find no minimal weak link of weight 7 with
$p<3$ suggests that $p$ continues to grow with $w$ for higher values,
and it may be that eqn. (\ref{eq_dsdt}) is valid for all $w$. 

By definition, a minimal weak link is on the boundary between being
`loosely connected' and `well connected'---deletion of an edge
would make it not weakly linked, addition of an edge would make it
either strong or not minimal. We find that
the number of edges grows roughly linearly with $w$. If one
constructs links simply by replacing edges by bridges, starting from W1,
then each time the number of vertices increases by two, the number of edges
increases by three, so one will obtain $(3w+1)/2$ edges. We have not
found a link with fewer edges than that. The average vertex degree for
a graph with this number of edges is $(3w+1)/(w+2)$, i.e. tending to
3 from below as $w \rightarrow \infty$. This construction leads to
a graph with many vertices of degree 2, and some of degree 4 or more.
It is interesting to ask whether minimal weak links of arbitrary $w$ have
to contain vertices of high degree. To begin an exploration of
this question, we offer the following observations.

A minimal strong link with any desired distance between the terminals
(where the distance between the terminals is defined as the shortest
walk between them on the graph) can be achieved with maximum vertex
degree 4, by using a `chain'. A `chain' SC$n$ is a sequence of mutually
supporting pairs having a structure of a generic form illustrated by SC3:

\begin{picture}(40,35)(30,-5)
\multiput(47.5,7.5)(130,0){2}{\rule{5pt}{5pt}}
\multiput(80,20)(0,-20){2}
{ 
  \multiput(0,0)(35,0){3}{\circle{5}}
  \multiput(2.5,0)(35,0){2}{\line(1,0){30}}
}
\multiput(82.2,19)(35,0){2}{\line(5,-3){30}}
\multiput(82.2,1)(35,0){2}{\line(5,3){30}}
\multiput(52.5,11)(100,-10){2}{\line(3,1){25}}
\multiput(52.5,9)(100,10){2}{\line(3,-1){25}}
\end{picture}\\
That is, SC1$\equiv$S2, SC2 is the link W5C with the pendant removed,
the diagram above shows SC3, and the extension to higher $n$ is obvious.
SC$n$ has maximal degree $d_{\rm max} = 2,3,4$ for $n=1$, $n=2$,
$n \ge 3$ respectively and the distance between terminals is $n+1$. By
adding a pendant one obtains a weak link of weight $2n+1$ with distance
between terminals equal to $n+2$. It is not hard to convince oneself
that this is the maximum distance possible with given weight. Two further
questions now suggest themselves: first, can there be weak links
of arbitrary weight and maximal degree less than 4? Secondly, what is the
largest possible distance between terminals for links with $d_{\rm max}=3$?
By replacing one or more edges between $d=2$ vertices in W5Z by a short
chain (SC1 or SC2), one can construct weak links with $d_{\rm max}=3$
of weight 7,9,11,13. However, we find from our exhaustive search
that there are no strong links of weight in the range $5 \le w \le 8$ which have
$d_{\rm max}=3$ and both terminals of degree 2. Also, for $w=9$
there is no weak link with $d_{\rm max}=3$ 
containing a pair of neighbouring $d=2$ vertices.
These facts together suggest (but do not prove) that if there is a weak link
of weight 15 having $d_{\rm max}=3$ then it is Shannon-irreducible. We 
conjecture that there is no
such weak link. Finally, we find from the search that, for weak links
having $w \le 9$ and $d_{\rm max}=3$, the maximum distance between terminals is 4.

\subsection{Ladders and an example in Hex}

Figure \ref{f.hex45}(a) shows an example of W5Z (to be precise, a link 
which is bridge-reducible to W5Z)
in the opening position of $5 \times 5$ Hex. This link suffices to prove
that the two dotted vertices are winning openings for \short. 

W7E and W7F both have natural extensions to
larger links which realise the same basic idea but with more `rungs' in the ladder.
Let us introduce the notation WLL$(r)$ and WLT$(r)$ for this general construction,
where $r$ is the number of rungs. We will use this notation to denote not
a single link but a set of links all having the same Shannon-reduction. 
WLL$(r)$ signifies a longitudinal ladder (one running from one terminal
to the other) and WLT$(r)$ signifies a transverse ladder (one running across
the neighbourhood of a terminal), in either case with further vertices
furnishing the forking ladder escape. Thus W7E is the Shannon-irreducible member of
the set WLL$(2)$ and W7F is the Shannon-irreducible member of the set WLT$(2)$.
WLL$(1)$ and WLT$(1)$ have the same structure; the Shannon-irreducible
member is in both cases W5Z. Figure \ref{f.hex45}(b) illustrates
the use of WLL$(3)$ to prove a simple property of $4 \times 4$ Hex.
The figure shows an example of
WLL$(3)$ occuring in $4 \times 4$ Hex (the pivot is strongly linked
to one terminal and weakly linked to the other; above it is a 3-rung ladder,
the last rung of which is a bridge (S2)). That the dotted vertex 
is a winning opening move for \short\ is more easily proved using
a smaller and simpler weak link in this Hex game. 
However, the link shown is useful for another
purpose. Clearly if \cut\ is the first player, an opening move at
the isolated vertex (i.e. the vertex not needed by this link) is a losing move.
This fact could not be discovered merely by finding Shannon-reducible weak
links in the opening position.


\begin{figure}
\myfig{0.2}{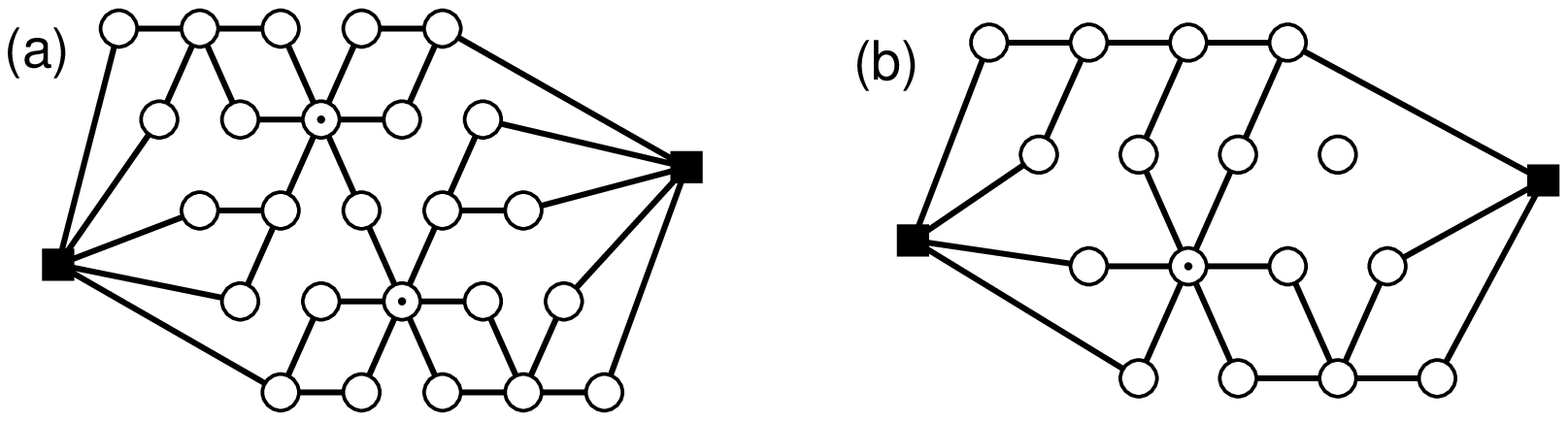}{f.hex45}
{Examples in Hex openings. (a) W5Z
in $5\times 5$ Hex; (b) WLL$(3)$ in $4\times 4$ Hex.
In each case, edges existing in the game but not needed for the given
link are not shown, in order to make the structure of the link easy to see.}
\end{figure}

\section{Conclusion}

We have presented constraints on graphs containing
minimal weak links, and used them to allow an exhaustive search for such links
on up to 11 vertices. We have elucidated the structure of these links by
proposing three simple reductions and examining the irreducible links.
The number of Shannon-irreducible links is small for $n \le 10$, which is
encouraging for the general task of solving Shannon games, but
it begins to grow rapidly at $n=11$, and the irreducible links
grow in complexity as $n$ increases, as one must expect
from the known intractability of the general problem. 
Further work in this direction could be devoted to tightening the constraints;
finding further constraints; proposing further useful reductions; and
suggesting strategies for detecting links when they are present as a sub-graph
(not necessarily vertex-induced) between arbitrary pairs of vertices in some
Shannon game. 



\section{Appendix}

{\em Proof of the case $w > 5$ in theorem \ref{th.dterm}.}

Suppose the converse, i.e. $p < 2$, in order to find a contradiction.
We already proved that $p \ge 1$ for $w>3$ so we infer that $p=1$.
That is, there is a single central vertex. Label it $c$.
We consider the following cases: (i) one of the terminals is pendant, (ii) one
of the terminals has neighbours all of degree 2, (iii) the rest---i.e. the borders
have at least two vertices each, and at least one vertex in each
border has degree greater than 2. 

Case (i). Let $s$ be the pendant terminal, and let
$v$ be its neighbour. This is the whole
of the first border, $B_s$. $v$ cannot be of degree 2 or
there is no weak link (\cut\ has a winning strategy). Therefore it is adjacent
to at least two other vertices. At most one of these can be in $B_t$,
or the link would contain W3. Therefore $\Gamma(v) = \{s, c, b\}$
where $b \in B_t$. Now consider the other members of $B_t$. There are at least
3 further members since $w>5$. Consider any pair of them. These cannot be of
degree 1, and cannot be adjacent to any other vertex in $B_t$ (or there would
be a dead edge), and cannot be adjacent to $v$ (or the link would contain W3).
Hence they are both adjacent to $c$. But then the link contains W5B, so it is
not minimal---we have a contradiction.

Case (ii). Let $s$ be a terminal whose neighbours all have degree 2.
All the members of $B_s$ are therefore threatened. This means
\cut\ has a winning strategy unless \short\ can remove two or more such
threats in a single move. This means there is a vertex not
in $\{s\} \cup B_s$ which is adjacent to at least two 
members of $B_s$. This vertex cannot be in $B_t$ or the link would contain
W3, therefore it is the vertex $c$. 
There must also be an edge between $c$ and some $u \in B_t$ or $c$
would be dead by theorem TDS.16.
So far we have proved that the link
contains at least the following vertices and edges:

\begin{picture}(40,35)(30,-5)
\multiput(47.5,7.5)(120,0){2}{\rule{5pt}{5pt}}
\put(52.5,11){\line(3,1){25}}
\put(52.5,9 ){\line(3,-1){25}}
\put(110,10){\circle{5}}
\multiput(110,10)(60,0){2}
{ \multiput(-30,9)(0,-18){2}{\circle{5}}
  \put(-2.5,1){\line(-3,1){25}}
  \put(-2.5,-1){\line(-3,-1){25}}
}
\put(112.5,11){\line(3,1){25}}
\put(48,0){$s$} \put(106,0){$c$} \put(173,0){$t$}
\put(145,20){$u$}
\end{picture}\\
Now consider another vertex $v \in B_t$, $v \ne u$.
If there is no such vertex then we have case (i) which we already proved.
$v$ cannot be adjacent any other member of $B_t$ (or there would be a dead edge)
and it cannot be adjacent to $c$ or the link would contain W5X. It cannot
be adjacent to the vertices of $B_s \cap \Gamma(c)$ since they have degree 2
by assumption. It follows that either $v$ is dead, or there are further
vertices in $B_s$ and it is adjacent to one of those. However, if it is
only adjacent to one of them then it is in a mutually threatening pair,
and if it is adjacent to more than one of them then the link contains W3. 
In either case the link is not minimal, and we have a contradiction.

Case (iii). Consider the number of edges between $B_s$ and $B_t$.
If there are none, then every vertex in $B_s$ must be adjacent to $s$ and $c$
but no other vertex (to avoid dead edges). Therefore we have case (ii) which
we already proved. If there is exactly one edge between $B_s$ and $B_t$ 
then we have that the link contains at least the following:

\begin{picture}(40,35)(30,-5)
\multiput(47.5,7.5)(120,0){2}{\rule{5pt}{5pt}}
\put(52.5,11){\line(3,1){25}}
\put(52.5,9 ){\line(3,-1){25}}
\put(110,10){\circle{5}}
\multiput(110,10)(60,0){2}
{ \multiput(-30,9)(0,-18){2}{\circle{5}}
  \put(-2.5,-1){\line(-3,-1){25}}
}
\multiput(112.5,9)(30,10){2}{\line(3,-1){25}}
\put(82.5,19){\line(1,0){55}}
\put(48,0){$s$} \put(106,0){$c$} \put(173,0){$t$}
\put(71,20){$u$} \put(145,20){$v$}
\end{picture}\\
But for the case under consideration, at least one vertex in each border
has degree greater than 2. These vertices can only be $u$ and $v$ if there are
no further edges between $B_s$ and $B_t$. Hence $u$ and $v$ must both be adjacent
to $c$, therefore the link contains W5X. It remains to consider the situation
where there is more than one edge between $B_s$ and $B_t$. We then have that the
link contains

\begin{picture}(40,35)(30,-5)
\multiput(47.5,7.5)(120,0){2}{\rule{5pt}{5pt}}
\put(52.5,11){\line(3,1){25}}
\put(52.5,9 ){\line(3,-1){25}}
\put(110,10){\circle{5}}
\multiput(80,20)(0,-21){2}
{ 
  \multiput(0,0)(60,0){2}{\circle{5}}
  \put(2.5,0){\line(1,0){55}}
}
\put(142.5,-1){\line(3,1){25}}
\multiput(82.5,19)(60,0){2}{\line(3,-1){25}}
\put(48,0){$s$} \put(114,6){$c$} \put(173,0){$t$}
\put(71,-6){$u$} \put(145,20){$v$} \put(145,-6){$w$}
\end{picture}\\
where we used that the link must not contain W3, and at least one of $B_s$
must have degree $>2$ (or we have case (ii)). Now, at least one of $v,w$
has a further edge or we have case (ii), and at least one of $u,w$
has a further edge, or they form a mutually threatening pair. Such an
edge cannot be within a border (or it would be dead) nor between borders (or
we have W3), therefore it must be to $c$. Therefore the link contains either
the edge $cw$ or both the edges $uc$ and $cv$. In either case the link contains
W5Z so is not minimal. This exhausts all possible contructions and so we
have a contradiction, which proves the theorem. \QED

{\em Proof that there is no minimal weak link of weight $w=4$.}

We will show that for a link of weight four, either \short\ has no
1st-player winnning strategy (i.e. it is not a weak link),
or the link contains W3 and hence is not minimal.
There can be at most one pendant, so the terminal degrees are either
$(1,2)$ or $(1,3)$ or $(2,2)$. In the case $(1,2)$ first suppose the
non-terminal vertices induce a complete graph (a clique). Now in order
to avoid containing W3, an edge from one border to the other must be deleted.
As soon as this is done, there is no weak link---\cut\ has a 2nd player
winning strategy. In the case $(1,3)$, in order to avoid a \cut\ win,
the neighbour to the pendant terminal must have degree $d\ge 3$. Therefore it
has at least two edges to the other border, so the link contains W3.
In the case $(2,2)$, to avoid a \cut\ win either both 
vertices in a border must be adjacent to one vertex in the other, in which
case the link contains W3, or at least one vertex in a border has degree $d\ge 3$,
in which case again the link contains W3. \QED

\begin{lemma}  \label{c_adjBt}
If a vertex $b \in B_s$ is adjacent to all the central vertices, then
every central vertex is adjacent to $B_t \setminus \Gamma(b)$.
\end{lemma}

Proof. Let $c$ be a central
vertex. By definition a central vertex is not adjacent to either terminal. 
$c$ must be adjacent to at least one vertex not in $\Gamma(b)$,
since otherwise the edge $bc$ is transverse and therefore dead. $c$ 
must also be adjacent to at least one vertex not in $B_s$, otherwise
it is surrounded by a terminal and thus dead by theorem TDS.16. 
It follows that must be an edge between $c$ and a member of
$B_t \setminus \Gamma(b)$. \QED

{\em Proof of the case $w>5$ in theorem \ref{th.demax}}.

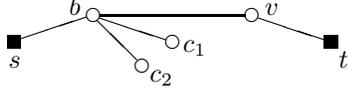
\begin{figure}
\begin{picture}(40,35)(30,-5)
\multiput(47.5,7.5)(120,0){2}{\rule{5pt}{5pt}}
\put(52.5,11){\line(3,1){25}}
\put(110,10){\circle{5}}
\multiput(80,20)(0,-21){1}
{ 
  \multiput(0,0)(60,0){2}{\circle{5}}
  \put(2.5,0){\line(1,0){55}}
}
\multiput(82.5,19)(60,0){2}{\line(3,-1){25}}
\put(81.5,18){\line(1,-1){15} \put(2,-17){\circle{5}} \put(5,-23){$c_2$}}
\put(48,0){$s$} \put(114,6){$c_1$} \put(173,0){$t$}
\put(71,20){$b$} \put(145,20){$v$}
\end{picture}
\caption{Subgraph used in proof of theorem \ref{th.demax}.}
\label{f.subg}
\end{figure}

We want to prove $d_b \le p+1$; we suppose the contrary, i.e. $d_b > p+1$,
and seek a contradiction. We already proved that $d_b \le p+2$ so we are
supposing $d_b = p+2$. Now, $b$ cannot be adjacent to other vertices in its
own border (or there would be dead edges) and it can be adjacent to at most
one vertex in the other border (or the link would contain W3). Hence it must be that
$b$ is adjacent to all the centre vertices, and to its terminal and to one vertex
in the other border. Also, we know from theorem \ref{th.dterm} that $p \ge 2$.
It follows that the link contains at least the vertices and edges shown
in figure \ref{f.subg}.

Now, by lemma \ref{c_adjBt}, there
must be edges from $c_1$ and $c_2$ to one or more members of $B_t \setminus \{v\}$.
If both are adjacent to the same vertex $w \in B_t$, $w \ne v$, then the
link contains W5A. If either one of them is adjacent to two vertices
in $B_t \setminus \{v\}$ then the link contains W5B. It follows that the link
contains the following structure, with no further edges between
$C \equiv \{c_1,c_2 \}$
and $B_t \setminus \{v\}$:

\begin{picture}(40,35)(30,-5)
\multiput(47.5,7.5)(120,0){2}{\rule{5pt}{5pt}}
\put(52.5,11){\line(3,1){25}}
\put(110,10){\circle{5}} \put(112,3){$c_1$}
\multiput(80,20)(0,-21){1}
{ 
  \multiput(0,0)(60,0){2}{\circle{5}}
  \put(2.5,0){\line(1,0){55}}
}
\multiput(82.5,19)(60,0){2}{\line(3,-1){25}}
\put(81.5,18){\line(1,-1){15} \put(2,-17){\circle{5}} \put(2,-25){$c_2$}}
\put(170,10)
{ \multiput(-30,0)(0,-9){2}{\circle{5}}
  \put(-2.5,-1){\line(-3,-1){25}}
  \put(-32.5,0){\line(-1,0){25} \put(5,0){\line(1,0){25}}}
  \put(-69,-9){\line(1,0){37}}
}
\put(48,0){$s$}  \put(173,0){$t$}
\put(71,20){$b$} \put(145,20){$v$}
\end{picture}\\
Suppose that $d_t > p+1$. Then $B_t$ contains further vertices. Each
of these is not
adjacent to any vertex in $\{b\} \cup C \cup B_t$ since otherwise
the link would contain one or both of W3, W5B, or there would be a dead edge.
Hence the remaining vertices in $B_t$ are only adjacent to $t$ and a member
of $B_s \setminus\{ b\}$, and each can be adjacent to at most one of the latter or the
link would contain W3. It follows that all these further members
of $B_t$ have degree 2.
We have that each further vertex $w \in B_t$ has degree 2 and
is connected to a further vertex $u \in B_s$ whose degree must exceed 2
or we would have a mutually threatening pair and the link would not be 
minimal. Therefore each such pair $uw$ is a short-cut. 
However, the further edge or
edges from $u$ (in addition to $us$ and $uw$) can only go to $C$ (to avoid
dead edges and W3). But this is ruled out by theorem \ref{shortcut_dead}---if
there were an edge $uc$ where $c \in C$, then the edge $bc$ would be a spectator.
Hence we have ruled out $d_t > p+1$. 

We now finish by showing that if $d_t= p+1$
the link cannot be weak. For, suppose $d_s > 1$, i.e. there are further
vertices in $B_s$. These may not be adjacent to each other or $b$ (to
avoid dead edges), nor to
$v$ (or the link contains W3), nor $\Gamma(C) \cap B_t$
(or the link contains W5Z), nor $t$ (to avoid W1). Therefore they
are only adjacent to $C$.
Now we allow \short\ to short all of them, i.e. all of $B_s \setminus \{ b\}$,
as a free move. The result is to introduce edges between $s$ and $C$,
and possibly within $C$,
and make $s$ a pendant. Therefore if \short\ has a winning strategy,
he must now open at $b$. \cut\ replies at $v$ and wins the game (all
the other vertices being now mutually threatening pairs after removal
of dead edges). We have a contradiction with the opening supposition. \QED

{\em Proof of theorem \ref{th.dmax}.}

Suppose there is a central vertex $c$ with $d_c = \max(d_s,d_t) + p$.
We will show that this is not possible for a minimal weak link
of weight $w>5$, with one exception. By definition, $c$ is not adjacent to either
terminal, and if it is adjacent to more than one vertex in one border
it can be adjacent to at most one vertex in the other border, or the
link would contain W5X. Hence if $d_c= \max(d_s,d_t) + p$ then
$c$ must be adjacent to
all the other central vertices, and to all the largest border, and
to one vertex in the other border. Therefore, the
link contains at least:

\begin{picture}(40,35)(30,-5)
\multiput(30.5,7.5)(137,0){2}{\rule{5pt}{5pt}}
\multiput(57,10)(-2,7){2}{\circle{5}}
\put(36,11){\line(3,1){17}}
\put(35,10){\line(3,0){19}}
\put(59,11){\line(5,2){19}}
\put(58,17){\line(6,1){19}}
\put(110,10){\circle{5}}
\multiput(80,20)(0,-21){1}
{ 
  \multiput(0,0)(60,0){2}{\circle{5}}
  \put(-4,-9){$c$}  \put(65,0){$v$}
  \put(2.5,0){\line(1,0){55}}
}
\multiput(82.5,19)(60,0){2}{\line(3,-1){25}}
\put(81.5,18){\line(1,-1){15} \put(2,-17){\circle{5}} \put(5,-23){$c_2$}}
\put(28,0){$s$} \put(114,6){$c_1$} \put(173,0){$t$}
\end{picture} \\
except that $c_2$ may be absent if $p=2$. First consider the case
$p > 2$ so we have this subgraph. This subgraph is the same as the
one shown in figure \ref{f.subg} after replacing the edge $sb$ by
a bridge (S2) and re-labelling $b \rightarrow c$. Therefore
the argument given above in the proof of the case $w>5$
of theorem \ref{th.demax}
can be used, as long as one keeps in mind that there may be further
edges between the vertices in this bridge and the rest of the graph.
Various minimal links appear in that argument. One and only
one of them can be constructed without spectator edges, namely
W5A. All other constructions either give a non-minimal link
or no link at all. Hence we find the link shown in figure
\ref{f.exlink} is an exception to the rule we are trying to prove,
and there are no other exceptions for $p > 2$. 

It remains to consider $p=2$. There is now just one central vertex in addition
to $c$; call it $c_1$. In order that $cc_1$ should not be transverse, we require
an edge from $c_1$ to a non-neighbour of $c$, which implies $B_t$ has at least
one vertex $w \ne v$, and $c_1$ is adjacent to $w$. If there are any more
vertices in the graph, then they must be border vertices and at least one of them
must be in $B_s$ since we assumed this is the largest border. If $d_s>2$
then let $b \in B_s$ be such a further vertex. It must be adjacent to $c$
(by the supposition that $c$ is adjacent to all the largest border). However
there is no edge $bv$ (or the link contains W5B), and there is no edge $bw$
(or the link contains W5Z expanded by the bridge, with $bc$ a spectator).
There must be an edge $c_1v$ or there is no weak link even after shorting all
of $B_s$, but in the presence of $c_1v$ there is no edge $bc_1$ or the link
contains W5X expanded by the bridge, with $bc$ a spectator again. We deduce
that $b$ is of degree 2. We may now apply the same reasoning to 
the other vertices in $B_s$,
showing that they also have degree 2. Hence we
have a group of three mutually supporting vertices, which is not possible
in a minimal link (lemma \ref{lem.support}). It follows that there is no
vertex $b$, so we have $d_s=2$ and therefore $d_t=2$ (since $d_s \ge d_t$
by assumption). It only remains to show that the link of weight 6 now
under consideration is either not weak or not minimal; this is easy to do.
\QED

\bibliographystyle{unsrt}
\bibliography{graphrefs}

\end{document}